\documentclass[reqno]{amsart}
\usepackage{xypic,amssymb}

\theoremstyle{plain}
\newtheorem{theorem}{Theorem} 
\newtheorem{lemma}[theorem]{Lemma}
\newtheorem{proposition}[theorem]{Proposition}

\theoremstyle{definition}
\newtheorem{defn}[theorem]{Definition}
\newtheorem{example}[theorem]{Example}

\theoremstyle{remark}
\newtheorem{remark}[theorem]{Remark}

\newcommand{\comment}[1]{}

\newcommand{\Pp}{{\mathbb P}}
\newcommand{\Zz}{{\mathbb Z}}
\newcommand{\Cc}{{\mathbb C}}

\newcommand{\EE}{{\mathcal E}}  
\newcommand{\FF}{{\mathcal F}}
\newcommand{\KK}{{\mathcal K}}
\newcommand{\LL}{{\mathcal L}}
\newcommand{\QQ}{{\mathcal Q}}
\newcommand{\WW}{{\mathcal W}}

\newcommand{\Tan}{T}           

\newcommand{\PV}{{\mathcal V}}               
\newcommand{\PS}{{\mathcal X}}               
\newcommand{\PSD}{{\mathcal X}^{d}}          
\newcommand{\PSM}[1]{{\mathcal X}^{M}(#1)}   
\newcommand{\PVD}{{\mathcal V}^{d}}          

\newcommand{\ideal}[1]{\left\langle #1 \right\rangle}
\newcommand{\inj}{\hookrightarrow}
\newcommand{\surj}{\twoheadrightarrow}

\newcommand{\Gr}{Gr}
\newcommand{\Fl}{{F\ell}}
\newcommand{\A}{{\mathcal A}}

\newcommand{\D}{{\mathcal D}}
\newcommand{\I}{{\mathcal I}}

\newcommand{\rank}{\mathop{\rm rank}\nolimits}
\newcommand{\Poin}{{\rm Poin}}
\newcommand{\PPD}[1]{P^{d}_{#1}(q)}       
\newcommand{\Tutte}{\mathbf{T}}
\newcommand{\defterm}[1]{{\it #1\/}}
\newcommand{\isom}{\cong}

\newcommand{\Pic}{\mathbf{P}}
\newcommand{\QPic}{\mathbf{Q}}
\newcommand{\sm}{\setminus}
\newcommand{\st}{~\mid~}       
\newcommand{\x}{\times}
\newcommand{\fld}{\mathbf{k}}
\newcommand{\Sch}{\mathfrak{S}}

\newcommand{\numv}{\mathbf{v}}
\newcommand{\nume}{\mathbf{e}}
\newcommand{\numc}{\mathbf{c}}
\newcommand{\numi}{\mathbf{i}}
\newcommand{\numl}{\mathbf{l}}

\newcommand{\dsum}{\oplus}

\begin{document}

\title{On the topology of graph picture spaces}
\author{Jeremy L. Martin}

\address{School of Mathematics\\
University of Minnesota\\
Minneapolis, MN 55455}
\email{martin@math.umn.edu}
\thanks{Supported by an NSF Postdoctoral Fellowship.}

\begin{abstract}
We study the space ${\mathcal X}^{d}(G)$ of
pictures of a graph $G$ in complex projective $d$-space.  The main result
is that the homology groups (with integer coefficients) of ${\mathcal
X}^{d}(G)$ are completely determined by the Tutte polynomial
of $G$.  One application is a criterion in terms of the Tutte polynomial
for independence in the {\it $d$-parallel matroids\/} studied in combinatorial
rigidity theory.  For certain special graphs called \defterm{orchards},
the picture space is smooth and has the structure of an
iterated projective bundle.  We give a Borel presentation of the cohomology
ring of the picture space of an orchard,
and use this presentation to develop an analogue of the classical Schubert calculus.
\end{abstract}

\maketitle


\section{Introduction} \label{intro-section}

The theory of \defterm{configuration varieties}, such as the Grassmannian,
flag and Schubert varieties, is marked by an interplay between different fields of 
mathematics, including algebraic geometry, topology and combinatorics.  In this paper,
we study a class of configuration varieties called \defterm{picture spaces of graphs},
a program initiated in~\cite{JLM1}.
As we will see, there is a close connection between the
combinatorial structure of a graph and the
topology and geometry of its picture space.

Let $G$ be a graph with vertices $V$ and edges $E$.
The \defterm{$d$-dimensional
picture space} $\PSD(G)$ is defined as the projective algebraic set
whose points are \defterm{pictures} of $G$ in complex projective
$d$-space~$\Pp^{d}=\Pp^{d}_{\Cc}$.  A picture~$\Pic$ consists of a
point $\Pic(v) \in \Pp^{d}$ for each vertex~$v$ of~$G$ and a line
$\Pic(e)$ for each edge~$e$, subject to the conditions $\Pic(v)
\in \Pic(e)$ whenever~$v$ is an endpoint of~$e$.

Two fundamental operations of graph theory are {\it deletion} and
{\it contraction}: given a graph $G$ and an edge $e$, we may delete
$e$ to form a graph $G-e$, or identify the endpoints of $e$
to form a graph $G/e$.  Many combinatorial invariants, such as the
number of spanning forests, the chromatic polynomial, etc., satisfy
a deletion-contraction recurrence; the most general and powerful
of these invariants is the {\it Tutte polynomial\/}
$\Tutte_G(x,y)$.  (For those unfamiliar
with the properties of the Tutte polynomial, we give a brief sketch
in Section~\ref{graph-subsection} below; a much more
comprehensive treatment may be found in~\cite{Brylawski}.)
In the context of the present study,
the graph-theoretic operations of
deletion and contraction correspond to canonical
morphisms~\eqref{epimorphism} and~\eqref{monomorphism}
between picture spaces.
This suggests that there is a connection between the geometry
or topology
of $\PSD(G)$ and the Tutte polynomial of $G$.

The main result of this paper (Theorem~\ref{main-thm} below)
characterizes the integral homology groups
of $\PSD(G)$ completely in terms of the Tutte polynomial.  Rather than attempting to
describe the homology directly, we apply topological machinery, such as
the Mayer-Vietoris sequence, to the morphisms arising from deletion
and contraction.
The result is a recurrence 
describing the
homology groups of $\PSD(G)$ in terms of those of
$\PSD(G-e)$ and $\PSD(G/e)$, where $e$ is any edge of $G$.  This recurrence may
in turn be phrased in terms of the Tutte polynomial, as follows.
%
%

\begin{theorem} \label{main-thm}
Let $G$ be a graph and $d \geq 2$ an integer.  Then
\begin{enumerate}
\item the picture space $\PSD(G)$ is path-connected and simply connected;
\item the homology groups $H_{i}(\PSD(G)) = H_{i}(\PSD(G);\:\Zz)$ are
free abelian for $i$ even and zero for $i$ odd; and
\item the ``compressed Poincar\'e series''
    $$
    \PPD{G} \ := \ \sum_{i} q^{i} \ \rank_\Zz \, H_{2i}(\PSD(G))
    $$
is a specialization of the Tutte polynomial $\Tutte_G(x,y)$, namely
    $$
    \PPD{G} \ = \
    ([d]_q-1)^{\numv(G)-\numc(G)} \ [d+1]_q^{\numc(G)} \ 
    \Tutte_G \left(\frac{\:[2]_q\:[d]_q}{[d]_q-1},\,[d]_q\right)
    $$
where $\numv(G)$ is the number of vertices of $G$, $\numc(G)$ 
is the number of connected components, and $[d]_q = (1-q^d)/(1-q)$.
\end{enumerate}
\end{theorem}

In Section~\ref{prelim-section}, we set forth some elementary facts and notation involving
graphs and picture spaces.  We assume some familiarity with basic graph theory,
for which~\cite{West} is an excellent reference (among many others).
For a more leisurely treatment
of the picture spaces of graphs, the reader is referred
to~\cite{JLM1}.

Section~\ref{main-proof-section} contains the proof of Theorem~\ref{main-thm}.
We continue by exploring some natural extensions of this main result.
In Section~\ref{manifold-section}, we consider the space $\PSM{G}$ of ``pictures''
of a graph $G$ on a complex manifold $M$.  With suitable conditions on $M$,
we may mimic the methods of Section~\ref{main-proof-section} to describe
the integral homology groups of $\PSM{G}$ in terms of the
Tutte polynomial of~$G$ and the dimension and Poincar\'e series of~$M$.

For a graph~$G$ without loops or parallel edges, the \defterm{picture variety}
$\PVD(G) \subset \PSD(G)$ is the (irreducible) algebraic variety
defined as the closure of the pictures $\Pic$ for which the points $\Pic(v)$ are distinct.
For more on this subject, seee~\cite{JLM1}.
The problem
of finding a combinatorial interpretation of the Poincar\'e series of
$\PVD(G)$ appears to be quite difficult. We work out a sample Poincar\'e series
calculation in Section~\ref{pvd-section}, and attempt to give some idea
of the obstacles that are likely to arise.

Section~\ref{parallel-section} describes an application of the main result to
the theory of \defterm{combinatorial rigidity}.  This subject (for which
an excellent reference is~\cite{GSS}) concerns questions such as the following.
Suppose that we are given a physical framework in $d$-dimensional space, built out of
``joints'' and ``bars'' corresponding to the vertices and edges of some graph $G$.
For our purposes, we suppose that the bars may
vary in length, but meet the joints at fixed angles.
How can one tell from the
combinatorial structure of $G$ whether the framework will hold its shape---whether
it is ``rigid'' or ``flexible''?  
We show in Section~\ref{parallel-section} that this information may be read off
from the Tutte polynomial specialization of Theorem~\ref{main-thm}.  In the language
of rigidity theory, this says that 
all of the information about the \defterm{$d$-parallel matroid} of a graph $G$
is contained in the Tutte polynomial $\Tutte_G(x,y)$.

In Section~\ref{orchard-section}, we study the 
multiplicative structure of the cohomology ring
$H^*(\PSD(G);\:\Zz)$ in the case that~$\PSD(G)$ is smooth.  This turns
out to be equivalent to the property that $G$ is an ``orchard'': that is,
every edge is either a loop or an isthmus.  For such a graph,
$\PSD(G)$ is an iterated projectivized vector bundle, so its cohomology ring has
a presentation in terms of Chern classes of line bundles (see~\cite{BT});
we give this presentation as Theorem~\ref{orchard-thm}.  This in turn
leads to a ``Schubert calculus of orchards'': that is, we can answer
certain questions involving the enumerative geometry of points and lines
in $\Pp^{d}_{\Cc}$ by means of polynomial calculations in $H^*(\PSD(G);\:\Zz)$.

The author thanks Wojciech Chach\'olski, Sandra Di~Rocco, 
and Victor Reiner for helpful discussions, and an
anonymous referee for providing numerous constructive suggestions
for improving this paper.


\section{Graphs and their Picture Spaces} \label{prelim-section}
\subsection{Graphs}  \label{graph-subsection}

We assume some familiarity with the basics of graph theory on the part of the
reader; a good general reference is~\cite{West}.

A \defterm{graph} is a pair $G = (V,E)$, where $V=V(G)$ is a
finite nonempty set of \defterm{vertices} and $E=E(G)$ is a set of
\defterm{edges}. Each edge $e$ has two vertices $v,w$, not
necessarily distinct, called its \defterm{endpoints}.  If $v=w$
then $e$ is called a \defterm{loop}.
When no ambiguity can arise, we sometimes denote an edge by its
endpoints, e.g., ``$e=vw$''. 
A \defterm{subgraph} $G'$ of $G$ is a graph with $V(G')
\subset V(G)$ and $E(G') \subset E(G)$.

The set of edges with~$v$ as an
endpoint is denoted~$E(v)$. The set $E(v) \cap E(w)$ is called a \defterm{parallel class}. A
graph is \defterm{simple} if all its of nonloop
parallel classes are singletons.  (In~\cite{JLM1}, there is the additional condition that
simple graphs contain no loops.  However, allowing loops does no harm to the results
of~\cite{JLM1} used here.)  An \defterm{underlying simple graph} of a graph $G$ is a
graph $G' = (V,E')$, where $E'$ consists of one member of each nonloop parallel class of $G$.

A graph $G$ is \defterm{connected} if for every $v,w \in V(G)$, there is a sequence of
vertices $v=v_0,v_1,\dots,v_r=w$ with $\{v_i,v_{i+1}\} \in E(G)$ for $0 \leq i < r$.  The
maximal connected subgraphs of $G$ are called its \defterm{connected components.}  We
write $\numv(G)$, $\nume(G)$, and $\numc(G)$ for, respectively, the number of vertices,
edges and connected components of $G$.

\begin{defn} \label{deletion-contraction-defn}
Let $e \in E(G)$.  The \defterm{deletion} $G-e$ is the graph
$(V,E \sm \{e\})$.  If $e$ is not a loop, the \defterm{contraction}
$G/e$ is obtained from $G-e$ by identifying the 
endpoints of $e$ with each other.
\end{defn}

For an example of contraction, see Example~\ref{messy-coincidence-locus-example}.  In general,
    \begin{align*}
    \nume(G-e) &\:=\: \nume(G)-1, & \nume(G/e) &\:=\: \nume(G)-1, \\
    \numv(G-e) &\:=\: \numv(G), & \numv(G/e) &\:=\: \numv(G)-1, \\
    \numc(G-e) &\:=\: \numc(G) ~\text{or}~ \numc(G)+1,  & \numc(G/e) &\:=\: \numc(G).
    \end{align*}
An edge $e$ is called an \defterm{isthmus} if $\numc(G-e) = \numc(G)+1$.

Some examples of graphs are as follows:
    \begin{itemize}
    \item the \defterm{complete graph $K_n$}, a simple graph with vertices 
        $\{1,\dots,n\}$, and edges $\{ij \st 1 \leq i<j \leq n\}$.
    \item the \defterm{empty graph} $N_n$, with $n$ vertices and no edges;
    \item the \defterm{loop graph} $L_1$, consisting of one vertex and a loop, and
    \item the \defterm{digon} $D_2$, with two vertices and two parallel nonloop edges.
    \end{itemize}

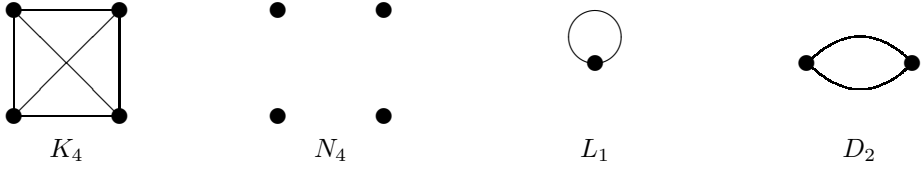
\begin{center}
\begin{picture}(340,70)
\put(20, 7){\makebox(0,0){$K_4$}}
\put( 0,20){\circle*{6}}
\put( 0,60){\circle*{6}}
\put(40,20){\circle*{6}}
\put(40,60){\circle*{6}}
\put( 0,20){\line(0,1){40}}
\put(40,20){\line(0,1){40}}
\put( 0,20){\line(1,0){40}}
\put( 0,60){\line(1,0){40}}
\put( 0,20){\line(1,1){40}}
\put( 0,60){\line(1,-1){40}}

\put(120, 7){\makebox(0,0){$N_4$}}
\put(100,20){\circle*{6}}
\put(100,60){\circle*{6}}
\put(140,20){\circle*{6}}
\put(140,60){\circle*{6}}

\put(220, 7){\makebox(0,0){$L_1$}}
\put(220,40){\circle*{6}}
\put(220,50){\circle{20}}

\put(320, 7){\makebox(0,0){$D_2$}}
\put(300,40){\circle*{6}}
\put(340,40){\circle*{6}}
\qbezier(300,40)(320,60)(340,40)
\qbezier(300,40)(320,20)(340,40)
\end{picture}
\end{center}

A fundamental isomorphism invariant of a graph $G=(V,E)$ is its
\defterm{Tutte polynomial} $\Tutte_G(x,y)$.
We describe here only a few
of the many properties of the Tutte polynomial; see the excellent survey by Brylawski and
Oxley~\cite{Brylawski} for more information.
For our present purposes, the following recursive definition of the Tutte polynomial
will be the most useful.

\begin{defn}
Let $G=(V,E)$ be a graph.  The \defterm{Tutte polynomial}
$\Tutte_G(x,y)$ is defined as follows.
If $\nume(G)=0$, then $\Tutte_G(x,y)=1$.  Otherwise, $\Tutte_G(x,y)$ 
is defined recursively as
    \begin{equation} \label{tutte-defn}
    \Tutte_G(x,y) \ = \ \begin{cases}
        x\cdot\Tutte_{G/e}(x,y) & \text{if}~e~\text{is an isthmus,} \\
        y\cdot\Tutte_{G-e}(x,y) & \text{if}~e~\text{is a loop,} \\
        \Tutte_{G-e}(x,y) + \Tutte_{G/e}(x,y) & \text{otherwise.}
    \end{cases}
    \end{equation}
for any $e \in E(G)$. (It is a standard fact, albeit not immediate from
the definition, that the choice of $e$ does not matter.)
\end{defn}

Many important graph isomorphism invariants satisfy deletion-contraction
recurrences of this form, and consequently may be obtained as specializations
of the Tutte
polynomial.  For instance, $\Tutte_G(1,1)$ equals the number of spanning
forests of $G$, while $\Tutte_G(2,2)=2^{\nume(G)}$.  In addition, one can
obtain more refined combinatorial data, such as the chromatic and flow
polynomials of $G$, by specializing the arguments $x$ and $y$ appropriately.
Again, the reader is referred to~\cite{Brylawski} for the full story.

There is an equivalent definition of the Tutte polynomial
as a certain generating function for the edge subsets $F \subset E$.
Define the \defterm{rank} of $F$, denoted $r(F)$,
as the cardinality
of a maximal acyclic subset of $F$; equivalently, $r(F)=\numv(G|_F)-\numc(G|_F)$,
where $G|_F$ is the edge-induced subgraph of $G$.
Then the Tutte polynomial may be defined in closed form as
the generating function
    \begin{equation} \label{tutte-gf}
    \Tutte_G(x,y) ~=~ \sum_{F \subset E} (x-1)^{r(E)-r(F)} (y-1)^{|F|-r(F)}
    \end{equation}
\cite[eq.~6.13]{Brylawski}; this formula will be useful when we study
the $d$-parallel matroid in Section~\ref{parallel-section}.


\subsection{Picture Spaces}  \label{picspace-subsection}

The main objects of our study are \defterm{picture spaces},
projective algebraic sets which
parameterize ``pictures'' of a graph in projective $d$-space
$\Pp^d$ over a field $\fld$.  In this paper, we shall be concerned
exclusively with the case $\fld=\Cc$; however, the picture space may be
defined over an arbitrary field.
The reader is referred to~\cite{JLM1}, especially Section~3, for
a more thorough discussion of the basic theory of picture spaces.

\begin{defn} \label{def-ps-one}
Let $G=(V,E)$ be a graph and $d \geq 2$ a positive integer.
A \defterm{picture} of $G$ in $\Pp^{d}$ is a tuple $\Pic$, consisting of a point 
$\Pic(v) \in \Pp^{d}$ for each $v \in V$ and a line $\Pic(e)$ in $\Pp^{d}$ for 
each~$e \in E$, such that $\Pic(v) \in \Pic(e)$ whenever $e \in E(v)$.
The set of all $d$-dimensional pictures is called the
\defterm{$d$-dimensional picture space} of $G$, denoted $\PSD(G)$.
\end{defn}

\begin{example}
A picture of $N_1$ is a point, so $\PSD(N_n) \isom \Pp^d$.  More generally,
a picture of the empty graph $N_n$ consists of an ordered
$n$-tuple of points in $\Pp^{d}$, so $\PSD(N_n) \isom (\Pp^d)^n$.

The data for a picture of the complete graph
$K_{2}$ consists of two points $\Pic(1),\Pic(2)$ and a line
$\Pic(12)$ containing both points. If $\Pic(1) \neq \Pic(2)$, then
$\Pic(e)$ is determined uniquely, but if $\Pic(1)=\Pic(2)$, then the
set of lines containing that point is isomorphic to $\Pp^{d-1}$.
In fact, $\PSD(K_{2})$ is the blowup of $\Pp^{d} \x \Pp^{d}$ along 
the diagonal $\{(p,p) \st p \in \Pp^{d}\}$ (see also~\cite[Example~3.6]{JLM1}).
This is a smooth, irreducible variety.
\end{example}

An easy consequence of the definition is that
if $G_1,\dots,G_r$ are the connected components of~$G$, then
$\PSD(G) \isom \PSD(G_1) \x \dots \x \PSD(G_r)$.  In particular, if
$\nume(G)=0$, then $\PSD(G) \isom (\Pp^d)^{\numv(G)}$.

Another elementary consequence is that if $e$ is a loop, then $\PSD(G)$ is a 
$\Pp^{d-1}$-bundle over $\PSD(G-e)$.  Indeed, a picture of $G$ may be regarded
as a picture of $G-e$, together with a line $\Pic(e)$ containing the point
$\Pic(v)$ (where $v$ is the single endpoint of $e$).  The fiber $\Pp^{d-1}$
corresponds to the space of lines through $\Pic(v)$ in $\Pp^d$.

This last observation may be generalized as follows.
Let $G=(V,E)$ be a graph and $G'=(V',E')$ a subgraph of $G$.
There is a natural epimorphism
    \begin{equation} \label{epimorphism}
    \PSD(G) \surj \PSD(G')
    \end{equation}
given by forgetting the picture data for vertices and edges not in $G'$.  Moreover,
if $G''=(V'',E'')$ is another 
subgraph of $G$, then the commutative diagram
    \begin{equation} \label{fiber-product}
    \xymatrix {
    & \PSD(G' \cup G'') \dlto\drto \\
    \PSD(G') \drto && \PSD(G'') \dlto \\
    & \PSD(G' \cap G'')
    } \end{equation}
is easily seen to be a fiber product square.  Here $G' \cup G''$ is the graph
with vertices $V(G') \cup V(G'')$ and edges $E(G') \cup E(G'')$, and $G' \cap G''$
is defined similarly.

Consider the Boolean algebra on $E$, where each subset $E' \subset E$ is associated with the
space $\PSD(V,E')$.  By (\ref{epimorphism}), there is an epimorphism $\PSD(V,E') \to
\PSD(V,E'')$ whenever $E'' \subset E'$.  Moreover, (\ref{fiber-product}) may be interpreted as
saying that the join of two spaces is the fiber product over their meet.  Accordingly,
$\PSD(G)$ is a fiber product of picture spaces of simple graphs---indeed, of graphs with one
edge each, which correspond to the atoms of the Boolean algebra.

\begin{example} \label{ps-example-one}
Consider the digon $D_2$, with vertices $1,2$ and parallel edges $e_1,e_2$.
A picture of $D_2$ consists of two points $\Pic(1),\Pic(2)$ and two lines
$\Pic(e_1),\Pic(e_2)$, such that $\Pic(i) \in \Pic(e_j)$ for $i,j \in \{1,2\}$.  
By the previous remarks, we may describe $\PSD(D_2)$ as a fiber product:
    $$\PSD(D_2) ~=~ \PSD(K_2) \;\underset{\PSD(N_2)}{\x}\; \PSD(K_2).$$
This space is neither smooth nor irreducible.  Its irreducible components are
    $$X_1 ~=~ \left\{ \Pic \in \PSD(D_2) \st \Pic(e_1)=\Pic(e_2) \right\},$$
which is isomorphic to $\PSD(K_2)$, and
    $$X_2 ~=~ \left\{ \Pic \in \PSD(D_2) \st \Pic(1)=\Pic(2) \right\},$$
which is a bundle over $\Pp^d$ with fiber $\Pp^{d-1} \x \Pp^{d-1}$.
The singular locus of $\PSD(D_2)$ is $X_1 \cap X_2$, which is isomorphic to a
$\Pp^{d-1}$-bundle over $\Pp^d$.
\end{example}

It will be useful to classify the pictures $\Pic$ of a graph
$G=(V,E)$ according to which points $\Pic(v)$ coincide.  Thus we are led to the
notion of a \defterm{cellule}.

\begin{defn} \label{define-cellule}
Let $\sim_{\A}$ be an equivalence relation on $V(G)$ with equivalence classes
$\A = \{A_{1}, \dots, A_s\}$.
The corresponding \defterm{cellule} in $\PSD(G)$ is defined as
    $$
    \PSD_{\A}(G) \ = \ \left\{ \Pic \in \PSD(G) \st \Pic(v)=\Pic(w) \;\iff\;
    v \sim_{\A} w \right\}.
    $$
A picture $\Pic$ is called \defterm{generic} if $\Pic(v) \neq \Pic(w)$ whenever
$v \neq w$.  Equivalently, $\Pic$ belongs to the \defterm{discrete cellule}
$\PSD_{\D}(G)$, where $\D$ is the equivalence relation
whose equivalence classes are all singletons. 
\end{defn}

Note that the cellules are
pairwise disjoint, and their union is~$\PSD(G)$.  Furthermore,
if $e$ is an edge whose endpoints
lie in different blocks of $\A$, then the points $\Pic(v)$ determine the line $\Pic(e)$
uniquely for each $\Pic \in \PSD_{\A}(G)$.  In this case we say that $e$ is
\defterm{constrained} with respect to $\A$.  Otherwise, varying $\Pic(e)$ while keeping the
other data of $\Pic$ fixed gives a family of pictures in $\PSD_{\A}(G)$;
this family is isomorphic to $\Pp^{d-1}$.
Therefore $\PSD_{\A}(G)$ has the structure of a fiber bundle, whose base is
$(\Pp^{d})^{|\A|}$ with diagonals deleted (that is, the discrete cellule of the
empty graph with vertices $V(G)$), and whose fiber is $(\Pp^{d-1})^{u(\A)}$, where
$u(\A)$ is the number of unconstrained edges.  In particular
    \begin{equation} \label{cellule-dimension}
    \dim\,\PSD_{\A}(G) \ = \ d|\A| \;+\; (d-1)\cdot u(\A)
    \end{equation}

In addition to (\ref{epimorphism}), there is a second canonical morphism between picture 
spaces, associated with any nonloop edge $e=vw$.
First, we define the \defterm{coincidence locus} $Z_{vw}(G)=Z_e(G)$ as
    \begin{equation} \label{coinc-locus}
    Z_{vw}(G) \ := \ \left\{ \Pic \in \PSD(G) \st \Pic(v)=\Pic(w) \right\}
    \ = \ \bigcup_{\A \; : \; v \sim_{\A} w} \PSD_{\A}(G).
    \end{equation}
Then there is a natural monomorphism
    \begin{equation} \label{monomorphism}
    \PSD(G/e) \inj \PSD(G-e)
    \end{equation}
whose image is the \defterm{coincidence locus} $Z_{vw}(G-e)$.

\begin{remark}
In light of (\ref{epimorphism}) and (\ref{monomorphism}), one may regard $\PSD$ as a
contravariant functor from the category of graphs to that of projective algebraic sets.  
Here a morphism $\phi:G \to G'$ of graphs is a pair of maps $V(G)\to V(G')$ and $E(G) \to
E(G')$ such that if $v$ is an endpoint of $e$, then $\phi(v)$ is an endpoint of $\phi(e)$.  
Thus $\PSD$ sends $\phi$ to a morphism
    $
    \phi^{\#} : \ \PSD(G') \to \PSD(G)
    $
defined by
    $$
    (\phi^{\#}\Pic)(a) \ = \ \Pic(\phi(a)),
    $$
where $\Pic \in \PSD(G')$ and $a$ is a vertex or an edge of $G$.  Furthermore,
if $G$ is a subgraph of $G'$, then $\phi^{\#}$ is an 
epimorphism, while if $G'$ is a quotient of $G$ (that is, it is obtained by a sequence of 
contractions), then $\phi^{\#}$ is a monomorphism.
\end{remark}


\section{Proof of the Main Theorem} \label{main-proof-section}

For the rest of the paper, we work over the ground field $\fld=\Cc$.
In this section, we show that the homology groups (with integer coefficients)
of $\PSD(G)$ are completely determined by the Tutte polynomial $\Tutte_G(x,y)$.
We have observed how the operations of deletion and contraction correspond
to the morphisms~\eqref{epimorphism} and~\eqref{monomorphism} of picture spaces.
In fact, these morphisms may be extended to a homotopy pushout
square~(\ref{pushout}), which induces a Mayer-Vietoris long exact
sequence of homology groups.  The Mayer-Vietoris sequence permits
us to write down a recursive formula for the Poincar\'e series of
$\PSD(G)$; this formula may in turn be expressed as a specialization
of the Tutte polynomial.  For the tools of topology that we use here,
an excellent reference is~\cite{Hatcher}.

We will work with the \defterm{$q$-analogues} $[q]_d$ of integers $d$, defined as
    $$[q]_d \ := \ (1-q^{d})/(1-q) \ = \ 1+q+\dots+q^{d-1}.$$

We shall see that for every graph $G$, the picture space $\PSD(G)$ has
only free abelian even-dimensional homology.  That is, the groups
$H_i(\PSD(G)) = H_i(\PSD(G);\:\Zz)$ 
are free abelian for all~$i$, and zero when $i$ is odd.  Accordingly, the
structure of the homology may be encoded conveniently by the
\defterm{compressed Poincar\'e series} whose coefficients are the even
Betti numbers:
    $$
    \PPD{G} \ := \ \sum_{i} q^{i} \ \rank_\Zz \, H_{2i}(\PSD(G)).
    $$

We first consider two simple cases.  If $\nume(G)=0$, then $\PSD(G) \isom
(\Pp^{d})^{\numv(G)}$, while if $\numv(G)=1$, then $\PSD(G)$ is a
$(\Pp^{d-1})^{\nume(G)}$-bundle over $\Pp^{d}$.  In both cases, $\PSD(G)$ is a simply
connected complex manifold with only free abelian even-dimensional homology.  Since the
compressed Poincar\'e series of $\Pp^{d}_{\Cc}$ is $[d+1]_q$, we have
    \begin{equation} \label{special-case}
    \PPD{G} \ = \ \begin{cases}
    [d+1]_q^{\numv(G)} & \quad \text{if}~ \nume(G)=0, \\
    [d]_q^{\nume(G)} \: [d+1]_q & \quad \text{if}~ \numv(G)=1.
    \end{cases}
    \end{equation}
In the case $\numv(G)=1$, we are using the fact that the Poincar\'e series of $\PSD(G)$ is
the same as that of $\Pp^{d} \x (\Pp^{d-1})^{\nume(G)}$ (see, e.g., Proposition~2.3
of~\cite{GaRe}).

Now let $G$ be an arbitrary graph and $e \in E(G)$ a nonloop edge.  As 
in~(\ref{epimorphism}) and~(\ref{monomorphism}), we have maps
    \begin{equation} \label{pre-pushout}
    \xymatrix{
    & \PSD(G) \ar@{->>}[d]^{\pi} \\
    \quad\PSD(G/e)\quad \ar@{^{(}->}[r]^{\iota} & \quad\PSD(G-e)\quad }
    \end{equation}
Since $\iota(\PSD(G/e)) \subset Z_e(G-e)$ and $\pi^{-1}(\PSD(G/e)) = Z_e(G)$, we may 
``complete the square'' to a commutative diagram
    \begin{equation} \label{pushout}
    \xymatrix{
    \quad Z_e(G) \quad \ar@{^{(}->}[r] \ar@{->>}[d] &
       \quad \PSD(G) \quad \ar@{->>}[d]^{\pi} \\
    \quad \PSD(G/e) \quad \ar@{^{(}->}[r]^{\iota} & \quad \PSD(G-e) \quad }
    \end{equation}

\begin{lemma} \label{hpo-lemma}
The map $\pi: Z_e(G) \to \PSD(G/e)$ is a $\Pp^{d-1}$-fibration, and the
diagram~(\ref{pushout}) is a homotopy pushout square.
\end{lemma}

\begin{proof}
Let $x,y$ be the endpoints of~$e$.
A picture in $Z_e(G)$ may be described by a picture of $G/e$, together with a line
in $\Pp^d$ through the point $\Pic(x)=\Pic(y)$.
Hence $\pi^{-1}(\Pic) \isom \Pp^{d-1}$.  In particular, $\pi$ induces an isomorphism between 
the fundamental groups of~$Z_e(G)$ and~$\PSD(G/e)$
(since complex projective $d$-space is simply 
connected).  Since the map $\iota$ is a monomorphism, the
diagram is a homotopy pushout.
\end{proof}

It follows from Lemma~\ref{hpo-lemma} that~(\ref{pushout}) induces a Mayer-Vietoris long 
exact sequence
    \begin{equation} \label{MV}
    \begin{array}{cccclcc}
    \dots &\to& H_{i}(Z_e(G)) &\to& H_{i}(\PSD(G/e)) \dsum H_{i}(\PSD(G)) &\to& 
    H_{i}(\PSD(G-e)) \\ \\
    &\to& H_{i-1}(Z_e(G)) &\to& \dots.
    \end{array}
    \end{equation}
This exact sequence allows us to compute
the homology groups recursively, in the same manner as the definition~\eqref{tutte-defn}
of the Tutte polynomial.  We can now prove the main theorem, which we
restate for convenience.

\vskip 0.1in

\noindent {\bf Theorem~\ref{main-thm}.}
Let $G$ be a graph and $d \geq 2$ an integer.  Then
\begin{enumerate}
\item the picture space $\PSD(G)$ is path-connected and simply connected;
\item the homology groups $H_{i}(\PSD(G))$ are free abelian for $i$
even and zero for $i$ odd; and
\item the compressed Poincar\'e series $\PPD{G}$ may be obtained from the
Tutte polynomial $\Tutte_{G}(x,y)$ by the formula
    \begin{equation} \label{poincare-formula}
    \PPD{G} \ = \
    ([d]_q-1)^{\numv(G)-\numc(G)} \ [d+1]_q^{\numc(G)} \ 
    \Tutte_G \left(\frac{\:[2]_q\:[d]_q}{[d]_q-1},\,[d]_q\right).
    \end{equation}
\end{enumerate}

\begin{proof}
We first show that $\PSD(G)$ is path-connected.  Let $\I$ be the
equivalence relation on $V$ in which all vertices are equivalent.
The corresponding
\defterm{indiscrete cellule} $\PSD_{\I}(G)$ consists of the pictures $\Pic$ with
$\Pic(v)=\Pic(w)$ for all $v,w \in V$.  Thus $\PSD_{\I}(G)$ is a
$(\Pp^{d-1})^{\nume(G)}$-bundle over $\Pp^d$; in particular it is path-connected.  On the
other hand, an arbitrary picture $\Pic$ can be deformed continuously into a picture in the
indiscrete cellule $\PSD_{\I}(G)$ as follows.  Choose a system of local affine coordinates
for $\Pp^d$ such that no $\Pic(v)$ lies on the hyperplane at infinity, rescale the
coordinates of all $\Pic(v)$ uniformly by a constant $\lambda$, and let $\lambda$ tend to
zero.  Since $\PSD_{\I}(G)$ is path-connected, so is $\PSD(G)$.

When $\numv(G)=1$ or $\nume(G)=0$, the formula for $\PPD{G}$ follows
from~(\ref{special-case}).  For the general case, we induct on $\numv(G)$ and $\nume(G)$.  In
particular, we may assume that the theorem holds for $\PSD(G-e)$ and $\PSD(G/e)$. Since
$Z_e(G)$ is a $\Pp^{d-1}$-bundle over $\PSD(G/e)$, it follows from Proposition~2.3
of~\cite{GaRe} that $Z_e(G)$ has only free abelian even-dimensional homology, and its
compressed Poincar\'e series is $[d]_q \; \PPD{G/e}$. Since~(\ref{pushout}) is a homotopy
pushout and the map $Z_e(G) \to \PSD(G/e)$ induces an isomorphism of fundamental groups,
the map $\pi$ also induces an isomorphism.  By induction, $\PSD(G)$ is simply connected.
Furthermore, the Mayer-Vietoris sequence~(\ref{MV}) gives $H_{i}(\PSD(G))=0$ for $i$ odd,
and splits into short exact sequences
    \begin{equation} \label{short-exact-sequence}
    0 \;\to\; H_{i}(Z) \;\to\; H_{i}(\PSD(G/e)) \dsum H_{i}(\PSD(G)) \;\to\; 
    H_{i}(\PSD(G-e)) \;\to\; 0
    \end{equation}
for $i$ even.  In particular, $H_{i}(\PSD(G))$ is free abelian for all~$i$, 

Recall the definition~(\ref{tutte-defn}) of the Tutte polynomial.  In order to
establish~(\ref{poincare-formula}), it suffices to show that
    \begin{subequations}
    \begin{align}
    \label{none} \PPD{G} &= [d+1]_q^{\numv(G)} && \text{if}~\nume(G)=0, \\
    \label{loop} \PPD{G} &= [d]_q \ \PPD{G-e} && \text{if}~e \in E~\text{is a loop,} \\
    \label{isth} \PPD{G} &= [2]_q \ [d]_q \ \PPD{G/e} && \text{if}~e~\text{is an isthmus,} \\
    \label{other} \PPD{G} &= \PPD{G-e} \:+\: ([d]_q-1)\PPD{G/e} && \text{otherwise.}
    \end{align}
    \end{subequations}
Indeed, (\ref{none}) is precisely~(\ref{special-case}), and (\ref{other}) follows
from~(\ref{short-exact-sequence}).  If $e$ is a loop, then $\PPD{G}$ is a $\Pp^{d-1}$-bundle 
over $\PPD{G/e}$, which implies~(\ref{loop}).  Now suppose that $e$ is
an isthmus with endpoints $v_1,v_2$.  It suffices to consider the case that $G$ is connected
and that $G-e$ has two connected components $G_1,G_2$, with $v_i \in V(G_i)$.  Then
    $$
    \PSD(G-e) \ \isom \ \PSD(G_1) \ \x \ \PSD(G_2)
    $$
and
    $$
    \PSD(G/e) \ \isom \ \PSD(G_1) \ \underset{\Pp^{d}}{\x} \ \PSD(G_2)
    $$
where the maps $\PSD(G_i) \to \Pp^{d}$ are given by $\Pic \mapsto \Pic(v_i)$ for $i=1,2$.
It follows that
    $$
    \PPD{G-e} \ = \ [d+1]_q \ \PPD{G/e}.
    $$
Substituting this into the recurrence defining $\PPD{G}$, we obtain
    $$
    \PPD{G} \ = \ ([d]_q-1) \ \PPD{G/e} \ + \ [d+1]_q \ \PPD{G/e}
    \ = \ [2]_q \ [d]_q \ \PPD{G/e}
    $$
yielding~(\ref{isth}) as desired.
\end{proof}


\section{Pictures on a Complex Manifold} \label{manifold-section}

We may generalize the preceding results by replacing $\Pp^d_\Cc$ with a complex manifold $M$
of dimension~$d$, provided that $M$ has only free abelian even-dimensional homology.  That is,
we study the space $\PSM{G}$ of ``$M$-pictures'' of $G$.  In this setting, we are faced
immediately with the problem of how to describe ``lines'' in~$M$ corresponding to edges
of~$G$. To resolve this question, we rely on the observations that for $e=vw \in E$, the data
$\Pic(e)$ is redundant if $\Pic(v) \neq \Pic(w)$, and if $\Pic(v)=\Pic(w)$, then $\Pic(e)$ may
be described as a tangent direction at $\Pic(v)$. With this in mind, we can formulate a
definition of $\PSM{G}$ which specializes to $\PSD(G)$ in the case that $M=\Pp^d$, and whose
Poincar\'e series obeys a recurrence akin to~(\ref{poincare-formula}).  We omit the proofs,
which are analogous to those of Theorem~\ref{main-thm}.

Denote by $\Tan M$ the tangent bundle of $M$, and by $\Tan_p M$ the tangent space of $M$ at a
point~$p$.  Note that $\Tan M$ is a rank-$d$ complex
vector bundle over~$M$, and its projectivization
$\Pp(\Tan M)$ (see Section~\ref{chern-subsection})
is a bundle over~$M$ with fiber~$\Pp^{d-1}$.

\begin{defn}
Let $G=(V,E)$ be a graph and $M$ a simply connected, compact complex manifold
of dimension~$d$ which has only free abelian even-dimensional homology.  The \defterm{space of 
$M$-pictures of $G$}, denoted $\PSM{G}$, is defined recursively as follows:
\begin{itemize}
\item If $E = \emptyset$, then $\PSM{G}$ is the direct product of $\numv(G)$ copies of~$M$.
\item If $G=L_1$, then $\PSM{G} = \Pp(\Tan M)$.
\item If $G=K_2$, then $\PSM{G}$ is the blowup of $M \x M$ along the diagonal.  (Note that 
blowups exist in the category of manifolds.)
\item Suppose that $E$ is the disjoint union of nonempty sets $E_1$ and $E_2$, and let 
$G_i=(V,E_i)$ for $i=1,2$.  Also, let $G_0=(V,\emptyset)$.  In this case, we define 
$\PSM{G}$ as a fiber product:
    $$
    \PSM{G} \ = \ \PSM{G_1} \; \underset{\PSM{G_0}}{\x} \; \PSM{G_2}.
    $$
\end{itemize}
\end{defn}

Let $e=vw \in E$.  As before, we define the coincidence locus $Z_e(\PSM{G})$ to be the set of
pictures $\Pic \in \PSM{G}$ for which $\Pic(v)=\Pic(w)$.  Then the analogue 
of~(\ref{pushout}) is the commutative diagram
    \begin{equation} \label{M-pushout}
    \xymatrix{ &
    \quad Z_e(\PSM{G}) \quad \ar@{^{(}->}[r] \ar@{->>}[d]^{\pi} &
    \quad \PSM{G} \quad \ar@{->>}[d] \\
    \quad \PSM{G/e} \ar@{{}*{\ \isom\ }{}}[r] &
    \quad Z_e(\PSM{G-e}) \ar@{^{(}->}[r] \quad & \quad \PSM{G-e} \quad
    } \end{equation}
which is a homotopy pushout square, with the map $\pi$ a $\Pp^{d-1}$-fibration. As in the
proof of Theorem~\ref{main-thm}, the Mayer-Vietoris sequence associated
with~(\ref{M-pushout}) implies (by induction) that $\PSM{G}$ has has only free abelian
even-dimensional homology, and leads to a recurrence for its Poincar\'e series.
This recurrence may in turn be translated into a formula in terms of the Tutte
polynomial of~$G$:
    \begin{multline}
    \sum_{i} q^{i} \ \rank_\Zz \ H_{2i}(\PSM{G}) \ = \\
    ([d]_q-1)^{\numv(G)-\numc(G)} \ \ 
    P(M)^{\numc(G)} \ \ 
    \Tutte_G \left(\frac{P(M)+[d]_q-1}{[d]_q-1},\,[d]_q\right)
    \end{multline}
where $P(M)=\Poin(M;\;q^{1/2})$.  One can easily check that this formula specializes
to that of Theorem~\ref{main-thm} in the case $M=\Pp^d$, $P(M)=[d+1]_q$.


\section{The Picture Variety} \label{pvd-section}

It is natural to ask if one can inductively calculate the Poincar\'e series of the picture
variety $\PVD(G)$, as we did above for the picture space $\PSD(G)$. However, it seems that the
only part of Theorem~\ref{main-thm} that can be extended to $\PVD(G)$ in general is the
simple argument for path-connectivity.  The problem is that the coincidence
locus of an edge in a picture variety is difficult to describe in general;
it is not always isomorphic to~$\PVD(G/e)$.  
The analogue of~\eqref{pushout} for picture varieties is the blowup diagram
    \begin{equation} \label{pushout-PV}
    \xymatrix{
    E \rto\dto^{p} & \PVD(G) \dto \\
    Z_e(G-e) \rto & \PVD(G-e) }
    \end{equation}
where $Z_e(G-e)$ now denotes the coincidence locus of $e$ in $\PVD(G-e)$
and $E$ is the exceptional divisor of the blowup (see Remark~3.5 of~\cite{JLM1}).
While~\eqref{pushout-PV} is a
homotopy pushout square, the coincidence locus $Z_e(G-e)$ may not be irreducible or even
equidimensional.  So the map $p$ in (\ref{pushout-PV}) may not be a
$\Pp^{d-1}$-fibration, and there appears to be no general way to calculate the Poincar\'e 
series of $\PVD(G)$ using a Mayer-Vietoris argument.

\begin{example} \label{pv-k4-example}
To give a sense of the difficulties involved, we sketch the calculation of the
Poincar\'e series of $\PVD(G)$
in the case $G=K_4$, $d=2$.  In a sense, this is the 
simplest case not covered by Theorem~\ref{main-thm}, since $K_4$ is the only simple graph 
with~$4$ or fewer vertices for which $\PV^2(G) \neq \PS^2(G)$.

Let $e$ be the edge joining vertices $1$ and $2$: then $\PV^2(G)$ is the blowup
of $\PV^2(G-e)$ along the coincidence locus of $e$.  This is the disjoint union of five
cellules, corresponding to the five equivalence relations $\A$ on $V(G)=\{1,2,3,4\}$
with $1,2$ equivalent.  Using the cellule dimension formula~\eqref{cellule-dimension},
one finds that one of these five cellules (the
indiscrete cellule) has codimension~$1$, three have codimension~$2$, and one has
codimension~$3$.  Every component of $Z$ must have codimension at most~$2$ (since $Z$ is
defined locally by two equations; see~\cite{JLM1}), so we can ignore the codimension-$3$
cellule.  On the other hand, blowing up a codimension-$1$ subvariety is a trivial operation,
so we can also ignore the indiscrete cellule.  Let $Z'$ be the union
of the three remaining cellules; then $\PV^2(G)$ is the
blowup of $\PV^2(G-e)$ along $Z'$, and the exceptional divisor $E'$ of the blowup
is a $\Pp^1$-bundle over $Z'$.  Moreover, replacing $Z,E$ with $Z',E'$ in~\eqref{pushout-PV},
we obtain a homotopy pushout square, and a Mayer-Vietoris sequence akin to~\eqref{MV}.

We can now
verify that $Z'$ is connected, simply connected, and has only free abelian even-dimensional
homology.  Moreover, the Poincar\'e series of $Z'$ may be calculated by applying
the Mayer-Vietoris sequence to its decomposition into three irreducible components.
By Proposition~2.3 of~\cite{GaRe},
$E'$ shares the first three properties, and its compressed Poincar\'e series
is $(1+q)$ times that of
$Z'$.  Finally, $\PV^2(G-e)=\PS^2(G-e)$ because $G-e$ is rigidity-independent
(see~\cite[Theorem 4.5]{JLM1}), so the Poincar\'e series of that space
is given by Theorem~\ref{main-thm}.
By an argument similar to the proof of Theorem~\ref{main-thm}, $\PV^2(G)$ has only free
abelian even-dimensional homology, and its compressed Poincar\'e series is
    $$
    (1+q)(1+q+q^{2})(q^{5}+13q^{4}+32q^{3}+24q^{2}+8q+1).
    $$
It is not clear what combinatorial significance this polynomial has, if any.
\end{example}

The technique of ignoring the codimension-1 components from the coincidence locus does allow us to
prove a general fact about picture varieties in the projective plane.

\begin{proposition} \label{simply-connected-two-prop}
Let $G$ be a simple graph.  Then $\PV^2(G)$ is simply connected.
\end{proposition}

\begin{proof}
Let $e \in E(G)$.
We have seen that in the blowup diagram~(\ref{pushout-PV}), every component of the coincidence 
locus $Z_e(G-e)$ has codimension~$1$ or $2$.  As in the preceding example, the 
codimension-$1$ components may safely be ignored, and the map
$p$ is a fibration with simply connected
fiber $\Pp^1$, hence induces an isomorphism of fundamental groups.  Since~(\ref{pushout-PV})
is a homotopy pushout, we also have an isomorphism between the fundamental groups of
$\PV^{2}(G)$ and $\PV^{2}(G-e)$, so we are done by induction on the number of edges.
\end{proof}

If $d>2$, then this argument does not go through because $Z'$ may not be equidimensional (the
codimensions of its components may vary between~$2$ and~$d$), in which case $p$ is not a
fibration.

A fundamental difficulty in extending our results to picture varieties is that
it is unclear how to define $\PVD(G)$ in the case that $G$ is not simple.  If we ``naively''
take $\PVD(G)$ to be the closure of the locus of generic pictures in $\PSD(G)$, then the
structure of parallel edges is lost. Indeed, if $e,e'$ are parallel, then $\Pic(e)=\Pic(e')$
for all generic pictures, hence for all pictures in $\PVD(G)$, which implies that $\PVD(G)$ is
isomorphic to the picture variety of its underlying simple graph.

A less trivial approach is to define $\PVD(G)$ as the locus of all pictures $\Pic \in \PSD(G)$
such that each ``underlying simple picture'' of $\Pic$ belongs to the picture variety of the
corresponding underlying simple graph.  This may be described as a fiber product, in the
spirit of~(\ref{fiber-product}).  Specifically, we define $\PVD(G)$ as the fiber product of
the picture varieties $\PVD(G_{i})$, where $G_{1},G_{2},\dots$ are the underlying simple
graphs of $G$. This definition preserves the structure of parallel edges; however, there still
exist coincidence loci which are hard to describe.

\begin{example} \label{messy-coincidence-locus-example}
Consider the graphs $G$ and $G/e$ given by the following figure.
    
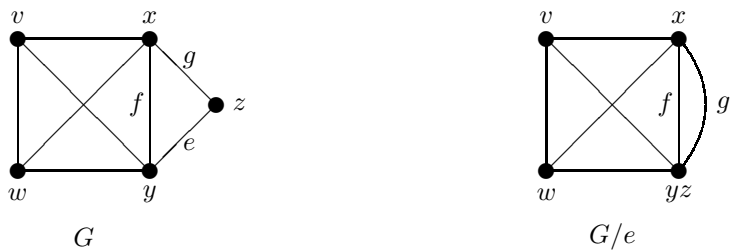
\begin{center}
    \begin{picture}(305,90)
    \put(45, 0){\makebox(0,0){$G$}}
    \put(20,25){\circle*{6}}
    \put(20,75){\circle*{6}}
    \put(70,25){\circle*{6}}
    \put(70,75){\circle*{6}}
    \put(95,50){\circle*{6}}
    \put(20,83){\makebox(0,0){$v$}}
    \put(70,83){\makebox(0,0){$x$}}
    \put(20,16){\makebox(0,0){$w$}}
    \put(70,16){\makebox(0,0){$y$}}
    \put(104,50){\makebox(0,0){$z$}}
    \put(20,25){\line(1,0){50}}
    \put(20,25){\line(1,1){50}}
    \put(20,25){\line(0,1){50}}
    \put(20,75){\line(1,0){50}}
    \put(70,25){\line(-1,1){50}}
    \put(70,25){\line(0,1){50}}
    \put(95,50){\line(-1,1){25}}
    \put(95,50){\line(-1,-1){25}}
    \put(85,35){\makebox(0,0){$e$}}
    \put(65,50){\makebox(0,0){$f$}}
    \put(85,67){\makebox(0,0){$g$}}
    \put(245, 0){\makebox(0,0){$G/e$}}
    \put(220,25){\circle*{6}}
    \put(220,75){\circle*{6}}
    \put(270,25){\circle*{6}}
    \put(270,75){\circle*{6}}
    \put(220,83){\makebox(0,0){$v$}}
    \put(270,83){\makebox(0,0){$x$}}
    \put(220,16){\makebox(0,0){$w$}}
    \put(270,16){\makebox(0,0){$yz$}}
    \put(220,25){\line(1,0){50}}
    \put(220,25){\line(1,1){50}}
    \put(220,25){\line(0,1){50}}
    \put(220,75){\line(1,0){50}}
    \put(270,25){\line(-1,1){50}}
    \put(270,25){\line(0,1){50}}
    \qbezier(270,25)(290,50)(270,75)
    \put(265,50){\makebox(0,0){$f$}}
    \put(287,50){\makebox(0,0){$g$}}
    \end{picture}
    \end{center}
Note that the nonparallel edges $f,g$ become parallel in $G/e$.  Let $G_1$ (resp.\ $G_2$) be 
the underlying simple graph of $G/e$ containing $f$ (resp.\ $g$), and let
$\A$ be the partition of $V(G)$ in which $\{y,z\}$ is the only nonsingleton block.
Since each $G_{i}$ is isomorphic to $K_4$,
Theorem~5.5 of~\cite{JLM1}) and our fiber product construction imply that the
picture variety $\PVD(G/e)$ is the vanishing locus in $\PSD(G/e)$ of two polynomials
$\tau(G_{1}),\tau(G_{2})$ (called in~\cite{JLM1} the \defterm{tree polynomials}
of $G_1$ and $G_2$ respectively).  Then
$\tau(G_{1})$ vanishes on $\PVD_{\A}(G-e)$ but $\tau(G_{2})$ does not, because $E(G) \sm
\{e,g\}$ is a rigidity circuit (see~\cite{JLM1}) while $E(G) \sm \{e,f\}$ is independent. Accordingly,
$\PVD_{\A}(G-e)$ is isomorphic neither to $\PSD(G/e)$ nor $\PVD(G/e)$.
\end{example}

We suspect that $\PVD(G)$ has in general only free abelian even-dimensional homology; it may
be possible to describe the centers of the iterated blowings-up sufficiently well to apply a
result such as Proposition~2.3 of~\cite{GaRe}.


\section{Applications to Parallel Independence} \label{parallel-section}

Let $\Pic$ be a $d$-dimensional picture of a graph $G=(V,E)$.  Consider a physical
model of $\Pic$ consisting of a ``bar'' for each edge $e$ and
a ``joint'' for each vertex $v$.  If $e$ has $v$ as an endpoint, then the corresponding
bar is attached to the corresponding joint.  The bars may cross, and their
lengths are allowed to vary, but we fix the angles at which the bars are attached
to the joints.  Thus, for example, a square framework may
be deformed to produce an arbitrary rectangle, but not any other rhombus.
Under what conditions on $G$ is such a model rigid?  That is, when is the model
determined up to scaling by specifying the attaching angles?
These and similar questions are the focus of \defterm{combinatorial
rigidity theory}; for more details, see, e.g.,~\cite{GSS} and~\cite{Whiteley}.

The graph $G$ (or, more properly, its edge set)
is said to be \defterm{$d$-parallel independent}
if for a generic picture in $\PSD(G)$, the directions of the lines
representing edges are mutually unconstrained.
For instance, $K_3$ is $2$-parallel
independent, because the slopes of two lines of a triangle in the plane
do not determine that of the third.  However, $K_3$ is $d$-parallel dependent
for $d\geq 3$, because the three sides of a triangle must be coplanar.
(The smallest simple graph which is $2$-parallel dependent is $K_4$.)

The condition of $d$-parallel independence is in fact a {\it matroid} independence
condition; for the reader not familiar with matroids, we remark here only that
it satisfies certain axioms which abstract the idea of linear independence
in a vector space.  Returning to rigidity for a moment,
a generic $d$-dimensional model of $G$, as described above, will hold its shape
if and only if $E(G)$ contains a $d$-parallel independent set of cardinality
$d\cdot\numv(G)-(d+1)$~\cite[Theorem 8.2.2]{Whiteley}.

The Poincar\'e series formula of Theorem~\ref{main-thm} may be applied
to give a criterion for $d$-parallel independence in terms of the Tutte
polynomial.  The central idea is that parallel independence can be determined from the
dimension and number of maximal-dimensional components of $\PSD(G)$, which can in turn be read
off from its Poincar\'e series.  Specifically,
let $X$ be a connected algebraic subset of $\Pp^n_\Cc$ of (complex) dimension $d$,
and let $c$ be the number of irreducible components of~$X$ of 
dimension~$d$.  By~\cite[Appendix A, Lemmas 2 and 4]{Ful97},
the leading term of the Poincar\'e series $\Poin(X;\:q)$ is~$cq^{2d}$.  We apply this fact
in the case that $X=\PSD(G)$.  (Note that the picture space is defined as an algebraic subset
of a product of Grassmannians, which may be identified in some complex projective space
by means of the Pl\"ucker and Segre embeddings (see, e.g.,~\cite[\S 9.1]{Ful97}.)

\begin{theorem} \label{d-parallel-thm}
Let $G=(V,E)$ be a graph and $d \geq 2$ an integer.  Then
$E$ is independent in the generic $d$-parallel matroid if and only if the polynomial
    $$
    ([d]_q-1)^{\numv(G)-\numc(G)} \ \Tutte_G \left( \frac{\:[2]_q\:[d]_q}{[d]_q-1},\,[d]_q\right)
    $$
is monic of degree $d(\numv(G)-\numc(G))$.
\end{theorem}

\begin{proof}
It suffices to consider the case that $G$ is simple,
since $L_1$ and $D_2$ are circuits in the $d$-parallel matroid for every~$d$
(this follows from Theorem 8.2.2 of~\cite{Whiteley}).  In this case,
Theorem~4.5 of~\cite{JLM1} tells us that $E$ is independent with respect to
the $2$-dimensional generic
rigidity matroid if and only if $\PV^2(G)=\PS^2(G)$.  The argument of~\cite{JLM1} can be
generalized as follows: for all $d$, $E$ is $d$-parallel independent if and only if
$\PVD(G)=\PSD(G)$. Recall that $\PVD(G)$ is an irreducible component of $\PSD(G)$ of dimension
$d \cdot \numv(G)$ and all other components have equal or greater dimension.
By the preceding remarks, $G$ is $d$-parallel independent if and only if the compressed
Poincar\'e series of $\PSD(G)$, as given by Theorem~\ref{main-thm}, is
monic of degree $d\cdot\numv(G)$.  By~\eqref{tutte-gf}, we have
    $$\Tutte_G\left(\frac{\:[2]_q\:[d]_q}{[d]_q-1},\,[d]_q\right)
    ~=~ \frac{f(q)}{([d]_q-1)^{r(E)}}
    ~=~ \frac{f(q)}{([d]_q-1)^{\numv(G)-\numc(G)}}$$
where $f(q)$ is a polynomial in $q$, and $r$ is the rank
function on subsets of $E$.  Therefore, dividing the
compressed Poincar\'e series by $[d+1]_q^{\numc(G)}$, which itself
is monic of degree $d\cdot\numc(G)$, produces a polynomial,
which is monic of degree $d(\numv(G)-\numc(G))$ if and only if
$G$ is $d$-parallel independent.
\end{proof}

One should note that this does not lead to an efficient algorithm for
computing the $d$-parallel behavior of an edge set.
The Tutte polynomial is exponentially hard to compute
(see~\cite{JVW}), and determining whether an edge set $F$
is $d$-parallel independent appears to be polynomial-time in $|F|$
(see, e.g.,~\cite[\S 4.10]{GSS}.



\section{Orchards} \label{orchard-section}

A graph $G=(V,E)$ is an \defterm{orchard} if every edge of $G$ is either a loop or an
isthmus.\footnote{To justify this terminology, $G$ is a forest with possibly some added loops,
which resemble fruit hanging from the trees.}
The orchards are precisely the graphs whose picture spaces are
smooth; indeed, the picture space of an orchard $G$ is an iterated bundle in
which each fiber is a complex projective space.  Thus we may write down an explicit ``Borel''
presentation for the integral cohomology ring $H^*(\PSD(G)) = H^*(\PSD(G);\:\Zz)$ in terms
of the Chern classes of natural line bundles on $\PSD(G)$.  The Borel presentation allows us
to formulate a ``Schubert calculus of orchards'', that is, a method of answering enumerative
geometry questions about~$\PSD(G)$ in terms of Schubert polynomial calculations.


\subsection{The Cohomology Ring} \label{chern-subsection}

Denote the number of isthmuses and loops of a graph $G$ by $\numi(G)$ and $\numl(G)$
respectively.  If $G$ is an orchard, then by~\eqref{tutte-gf} its Tutte polynomial is
$\Tutte_G(x,y) \ = \ x^{\numi(G)} y^{\numl(G)}$,
so by Theorem~\ref{main-thm} the compressed Poincar\'e series of $\PSD(G)$ is
    \begin{equation}
    \PPD{G} \ = \ [d+1]_q^{\numc(G)} \ [2]_q^{\numi(G)} \ [d]_q^{\nume(G)}.
    \end{equation}
This polynomial is palindromic, suggesting that the picture space of an orchard is smooth
(by Poincar\'e duality).  In fact, more is true.

\begin{proposition} \label{smooth-orchard-prop}
Let $G=(V,E)$ be a graph and $d \geq 2$.  The picture space $\PSD(G)$ is smooth 
if and only if $G$ is an orchard.
\end{proposition}

\begin{proof}
Suppose that $G$ is an orchard.  If $\nume(G)=0$, then $\PSD(G)=(\Pp^d)^{\numv(G)}$
is smooth.  Otherwise, let $e \in E$.  We shall shortly prove in 
Theorem~\ref{orchard-thm}
that $\PSD(G)$ is the projectivization of a complex vector bundle with
base $\PSD(G-e)$ (if $e$ is a loop) or $\PSD(G/e)$ (if $e$ is an isthmus).  It follows by
induction on $\nume(G)$ that $\PSD(G)$ is smooth.

Now suppose that $G$ is not an orchard.  It suffices to consider the case that $G$ is
connected and has no loops (since $\PSD(G)$ is the Cartesian product of the picture spaces of
the connected components of~$G$, and loops correspond to $\Pp^{d-1}$-bundles).  In this case
$\nume(G) \geq \numv(G)$.  Let $\Pic$ be a generic picture, so that $\PSD(G)$ looks locally
like $(\Pp^d)^{\numv(G)}$ near $\Pic$.  Therefore
    \begin{equation} \label{tanspace1}
    \dim \Tan_{\Pic}(\PSD(G)) \ = \ d \cdot \numv(G).
    \end{equation}  
At the other extreme, let $\QPic$ be a picture such that
$\QPic(v)=\QPic(w)=p$ for all $v,w \in V$ and $\QPic(e)=\QPic(f)=\ell$ for all $e,f \in E$.  
Then we can deform $\QPic$ into another picture of $G$ in the following ways:
    \begin{itemize}
    \item for any $v \in V$, move the point $\QPic(v)$ along the line $\ell$;
    \item for any $e \in E$, rotate $\QPic(e)$ through the space of lines containing $p$; or
    \item move the line in a direction ``orthogonal'' to itself (producing another totally 
        degenerate picture).
    \end{itemize}
(Here ``orthogonal'' really means ``orthogonal with respect to local affine coordinates in 
which the hyperplane at infinity does not contain $p$''.)  It follows that
    \begin{equation} \label{tanspace2}
    \begin{aligned}
    \dim \Tan_{\QPic}(\PSD(G)) \ &\geq \ \numv(G) \:+\: (d-1) \cdot \nume(G) \:+\: (d-1) \\
    &\geq \ d \cdot \numv(G) + d-1 \\
    &> \ d \cdot \numv(G),
    \end{aligned}
    \end{equation}
which together with (\ref{tanspace1}) implies that $\PSD(G)$ is not smooth.
\end{proof}

\begin{remark}
Let $G$ be the ``acetylene'' graph
\begin{picture}(30,20)
\put(5,5){\circle*{4}}
\put(25,5){\circle*{4}}
\put(5,5){\line(1,0){20}}
\qbezier(5,5)(15,15)(25,5)
\qbezier(5,5)(15,-5)(25,5)
\end{picture}.
This is not an orchard, so $\PS^2(G)$ is not smooth.  However, its Poincar\'e series
is palindromic: the formula of Theorem~\ref{main-thm} yields
$P^2_G(q) = 1+5q+9q^2+9q^3+5q^4+q^5$.
\end{remark}

In order to give a presentation of the cohomology ring $H^*(\PSD(G))$, we will need several
facts about vector bundles over complex manifolds.  We summarize these briefly. (For more
details, see chapter~IV of~\cite{BT}, especially pp.~269--271.)

Let $M$ be a complex manifold and $\EE$ a complex vector bundle of rank~$d$ over~$M$. The
\defterm{projectivization} of $\EE$ is the fiber bundle $\Pp(\EE) \overset{\pi}\to M$ whose
fiber at a point $m \in M$ is $\Pp(\EE)_m = \Pp(\EE_m)$, that is, the space of lines through
the origin in the fiber of $\EE$ at $m$.  Thus $\pi^{-1}\EE$ is a rank-$d$ vector bundle
over~$\Pp(\EE)$.  The \defterm{tautological subbundle} $\LL$ is the line bundle on $\Pp(\EE)$
defined fiberwise by $\LL_p = p$ (regarding $p$ as a line in $\EE_{\pi(p)}$).  Thus we have a
commutative diagram
    \begin{equation} \label{chern-setup}
    \xymatrix{ \LL \rto\drto & \pi^{-1}\EE \rto\dto & \EE \dto \\
    & \Pp(\EE) \ar@{->}[r]^{\pi} & M } \qquad
    \end{equation}
and with this setup, one has
    \begin{equation} \label{bundle-coho}
    H^*(\Pp(\EE)) \ \isom \ H^*(M)[x] \, / \,
    \ideal{ x^{d} + c_{1}(\EE)x^{d-1} + \dots + c_d(\EE) }
    \end{equation}
where $c_i(\EE)$ denotes the $i$th Chern class of~$\EE$, and
$x = c_{1}(\LL^*)$.  The Chern classes satisfy the following properties 
(see~\cite[p.~271]{BT}).  First,
if $\EE$ is a trivial bundle then $c_i(\EE)=0$ for $i>0$.
Second, if $\pi: X \to Y$ is a map of spaces and $\EE$ is a vector bundle on $Y$, then
    \begin{equation} \label{chern-pullback}
    c_i(\pi^{-1}\EE) = \pi^*(c_i(\EE)).
    \end{equation}
Third, there is the \defterm{Whitney product formula}: if
$    0 \,\to\, \EE' \,\to\, \EE \,\to\, \EE'' \,\to\, 0$
is an exact sequence of vector bundles on~$X$, then
    \begin{equation} \label{whitney}
    c(\EE) \ = \ c(\EE')\,c(\EE'')
    \end{equation}
where $c(\EE) = 1+c_{1}(\EE)+c_{2}(\EE)+\dots$ is the \defterm{total Chern class}.  In 
particular, if $\LL$ is a line bundle with dual bundle $\LL^*$, then $\LL \otimes \LL^*$
is a trivial bundle, so
    \begin{equation} \label{chern-dual}
    c_i(\LL^*) = -c_i(\LL).
    \end{equation}

Given an orchard~$G$, we may ``prune'' $G$ to obtain a smaller orchard~$G'$ such that
$\PSD(G)$ is the projectivization of a vector bundle $\EE$ on $\PSD(G')$.  The tautological
line bundle of $\EE$ may be expressed fiberwise in terms of the data $\Pic(v)$ and $\Pic(e)$,
allowing us to describe the cohomology ring of $\PSD(G)$ inductively by means
of~(\ref{bundle-coho}).

Denote by $\Gr(r,\Cc^{n})$ the Grassmannian variety of $r$-dimensional vector subspaces of
$\Cc^{n}$.  Thus $\Pic(v)$ and $\Pic(e)$ may be regarded as elements of $\Gr(1,\Cc^{d+1})$ and
$\Gr(2,\Cc^{d+1})$, respectively.  All the vector bundles we shall consider on $\PSD(G)$ are
subbundles of the trivial bundle $\WW = \Cc^{d+1} \x \PSD(G)$.  For each $v \in V$, there is a
line bundle $\LL_v \subset \WW$ with fiber
    \begin{equation} \label{lfiber}
    (\LL_{v})_{\Pic} \ = \ \Pic(v)
    \end{equation}
and for each edge $e \in E(v)$, there is a line bundle $\KK_{e,v}$ with fiber
    \begin{equation} \label{kfiber}
    (\KK_{e,v})_{\Pic} \ = \ \Pic(e) / \Pic(v).
    \end{equation}
The Chern classes of these bundles generate the cohomology ring of $\PSD(G)$,
as we now prove.

\begin{theorem} \label{orchard-thm}
Let $G=(V,E)$ be an orchard.
Then $H^*(\PSD(G);\:\Zz) \isom R/I$, where
    $$
    R \ = \ \Zz\big[x_v,y_{e,v} \ : \ v \in V, \: e \in E(v) \big]
    $$
and $I$ is the ideal
    $$\left\langle
    \begin{array}{ll}
    x_{v}^{d+1} & \qquad \text{for } v \in V, \\
    h_d(x_{v},\,y_{e,v}) & \qquad \text{for } v \in V, \; e \in E(v), \\
    x_{v}-x_{w}+y_{e,v}-y_{e,w}, \ 
    x_{v}y_{e,v}-x_{w}y_{e,w} & \qquad \text{for } e = vw
    \end{array}
    \right\rangle$$
where $h_d(x,y) = x^d+x^{d-1}y+\dots+xy^{d-1}+y^d$.
\end{theorem}

\begin{proof}
It suffices to consider the case that~$G$ is connected,
because the picture space of a graph is the product of the picture spaces of its 
connected components.

We induct on $\nume(G)$.  If $E = \emptyset$, then
    \begin{eqnarray*}
    H^*(\PSD(G)) &\isom& \Zz[x_v\ : \ v \in V] \ / \ 
        \ideal{ x_{v}^{d+1} \ : \ v \in V } \\
    &\isom& \bigotimes_{v \in V} \:\Zz[x_v] \ / \ \ideal{ x_{v}^{d+1} }
    \end{eqnarray*}
(the tensor product over~$\Zz$), because $\PSD(G) 
\isom (\Pp^d)^{\numv(G)}$.  Furthermore, $x_{v}$ is the dual Chern 
class of the line bundle with fiber $\Pic(v)$ (that is, the pullback to $\PSD(G)$ of
the tautological line bundle on the copy of $\Pp^d$ indexed by $v$).

Now, suppose that $\nume(G)>0$ and that the theorem holds for all orchards with fewer
edges than~$G$.  First, consider the case that there is at least one loop $e$, incident to a
vertex $v$.  Let $G'=G-e$.  Consider the vector bundles $\LL_v$ and $\WW$ on $\PSD(G')$, of
ranks~$1$ and~$d+1$ respectively.  A picture of $G$ may be specified by giving a picture
$\Pic$ of $G'$, together with a plane $\Pic(e) \in \Gr(2,\Cc^{d+1})$ such that $\Pic(v) 
\subset \Pic(e)$.  Equivalently, $\PSD(G) = \Pp(\QQ)$, where $\QQ$ is the
quotient bundle $\WW/\LL_v$.  Let $\KK_{e,v}$ be the tautological line bundle associated to
$\pi^{-1}\QQ$, so we have a diagram
    $$
    \xymatrix{
    \KK_{e,v}\rto\drto & \pi^{-1}\QQ\dto & \LL_{v}\rto\drto & \WW\rto\dto & \QQ\dlto \\
    & \PSD(G)=\Pp(\QQ) \ar@{->}[rr]^{\pi} && \PSD(G')}
    $$
Note that the map $\pi$ is a $\Pp^{d-1}$-fibration.
Let $y_{e,v} = c_{1}(\KK_{e,v}^*)$.  By~\eqref{bundle-coho},
    $$
    H^*(\PSD(G)) \ = \ H^*(\PSD(G'))[y_{e,v}] \ / \
    \ideal{ y_{e,v}^{d} + c_{1}(\QQ) y_{e,v}^{d-1} + \dots + c_d(\QQ) }.
    $$
By the Whitney formula, $c(\LL_v)c(\QQ) = c(\WW) = 1$.
Also, $c(\LL_v)=1-x_{v}$, so
    \begin{equation*}
    H^*(\PSD(G)) \ \isom \ H^*(\PSD(G'))[y_{e,v}] \ / \
        \ideal{ x_{v}^{d+1},\,x_{v}^{d}+x_{v}^{d-1}y_{e,v}+\dots+y_{e,v}^{d} }
    \end{equation*}
as desired.

Now, suppose that $G$ has a vertex $v$ with $E(v)=\{e\}$. Since $G$ is
connected, $e$ is not a loop. Let $w$ be the other endpoint of $e$, and let $G'$ be the
orchard obtained from $G$ by deleting $v$ and $e$ and adding a loop $\tilde{e}$ incident
to $w$. Forgetting the coordinates of $\Pic(v)$ and setting $\Pic(\tilde{e})=\Pic(e)$ gives an
epimorphism $\PSD(G) \to \PSD(G')$.  Moreover, a picture of $G$ may be specified by giving a
picture $\Pic$ of $G'$, together with a line $\Pic(v) \subset \Pic(e)$.  Hence $\PSD(G)$ is
the projectivization of the rank-$2$ bundle $\FF \to \PSD(G')$ whose fiber is $\FF_{\Pic} =
\Pic(e)$, and whose tautological subbundle is $\LL_v$.  We thus have a diagram analogous
to~(\ref{chern-setup}):
    $$
    \xymatrix{
    \LL_{v}\rto\drto & \pi^{-1}\FF\rto\dto & \FF\dto \\
    & \PSD(G)=\Pp(\FF) \ar@{->}[r]^{\pi} & \PSD(G') }
    $$

Let $x_{v} = c_1(\LL_{v}^*)$.  Then $H^*(\PSD(G))$ is generated by $x_v$ as an algebra
over $H^*(\PSD(G'))$.  Let $\KK_{e,v} = \FF/\LL_{v}$ and $y_{e,v} = c_1(\KK_{e,v}^*)$.  
By the Whitney formula,
$c(\FF) = c(\LL_{w})c(\KK_{e,w}) = c(\LL_{v})c(\KK_{e,v})$,
that is,
    $$
    (1-x_{v})(1-y_{e,w}) \ = \ (1-x_{w})(1-y_{e,v}).
    $$
Extracting the homogeneous parts of this equation, we find that
$x_{v}-x_{w}+y_{e,v}-y_{e,w}$ and $x_{v}y_{e,v}-x_{w}y_{e,w}$
are zero in $H^*(\PSD(G))$.  Eliminating $y_{e,v}$ from these equations recovers the
presentation~(\ref{bundle-coho}) of the cohomology ring, so it remains only to check the 
equations
    \begin{equation} \label{other-stuff}
    x_{v}^{d+1}=0, \qquad h_d(x_{v},y_{e,v})=0
    \end{equation}
in $H^*(\PSD(G))$.  Consider the loop graph $L_1$ with vertex $v'$ and edge $e'$.
Setting $\Pic(v')=\Pic(v)$ and $\Pic(e')=\Pic(e)$ gives an epimorphism $\pi: \PSD(G) \to 
\PSD(L_1)$.  From the first part of the proof, we have a presentation
    $$
    H^*(\PSD(L_1)) \ \isom \ \Zz[x,y] \:/ \ideal{ x^{d+1},\,h_d(x,y) }.
    $$
Here $x=c_1(\LL^*)$, where $\LL$ is the line bundle with fiber $\Pic(\tilde{v})$,
and $y=c_1(\KK^*)$, where $\KK$ is the line bundle with fiber
$\Pic(\tilde{e})/\Pic(\tilde{v})$.
Then $\LL_v = \pi^{-1}\LL$ and $\KK_{e,v} = \pi^{-1}\KK$, so the desired
equations~\eqref{other-stuff} follow from~\eqref{chern-pullback}.
\end{proof}

A more concise presentation of the cohomology ring can be obtained
by using the linear relations
$x_{v}-x_{w}+y_{e,v}-y_{e,w}$ to eliminate variables.  The most symmetric way to do this is to
introduce a new variable $z_e = x_{v}+y_{e,v} = x_{w}+y_{e,w}$
and then eliminate $y_{e,v}$ and $y_{e,w}$.  The resulting presentation is
as follows:
    \begin{multline} \label{concise-coho}
    H^*(\PSD(G)) \ \isom \ \Zz\left[x_{v},z_e \ : \ v \in V,e \in E\right] \ {\big /} \\
    \left\langle \begin{array}{rl}
        x_{v}^{d+1} & \quad:\quad v \in V, \\
        h_d(x_{v},\,z_{e}-x_{v}) & \quad:\quad e \in E(v), \\
        (x_{v}-x_{w})(z_{e}-x_{v}-x_{w}) & \quad:\quad e = vw
        \end{array} \right\rangle.
    \end{multline}


\subsection{Orchard Schubert Calculus} \label{schubert-subsection}

In this section, we give some examples of how Theorem~\ref{orchard-thm} may be
used to answer enumerative geometry questions,
in the spirit of the classical Schubert calculus on
Grassmannian and flag varieties.  Briefly, we can find the number of pictures
of a given orchard $G$ in $\Pp^d_{\Cc}$ satisfying certain incidence conditions,
by means of polynomial calculations in $H^*(\PSD(G))$.

We begin with a brief summary of the classical
theory; for more details, see~\cite{Ful97} or~\cite{KL72}.
Let $S_n$ be the symmetric group on $n$ letters, and let $\Fl(n)$
be the complete flag variety
    $$
    \Fl(n) \ = \ \left\{ F_{\bullet} = F_0 \subset \cdots \subset F_{n} \st F_i \in 
    \Gr(i,\Cc^n) \right\},
    $$
a complex manifold of dimension $n(n-1)/2$.  Fix a flag $F_{\bullet} \in \Fl(n)$.
For each $w \in S_n$, there is a corresponding \defterm{Schubert cell}, consisting of flags 
``in position $w$ with respect to $F_{\bullet}$'', that is,
    \begin{equation} \label{schubert-decomp}
    \Omega_\sigma^{\circ} \ = \ \left\{ E_{\bullet} \in \Fl(n) \st
    \dim(E_p \cap F_q) = \#\{i \leq p \st w(i) \geq n+1-q\} \right\}.
    \end{equation}
One has $\Omega_\sigma^{\circ} \isom \Cc^{n(n-1)/2-\ell(w)}$,
where $\ell(w) = |\{i<j \st w(i)>w(j)\}|$. The flag variety is the disjoint union of the
Schubert cells.  The \defterm{Schubert variety} $\Omega_{\sigma}$, defined as
the closure of $\Omega_{\sigma}^{\circ}$, is a union of Schubert cells:
    \begin{equation} \label{cell-closure}
    \Omega_{\sigma} \ = \ \overline{\Omega_\sigma^{\circ}} \ = \ \bigcup_{\rho \geq 
    \sigma} X_\rho^{\circ}
    \end{equation}
where $\geq$ is a certain partial order on $S_n$, the \defterm{strong Bruhat order}.
The cohomology classes $[X_\sigma] \in H^{2\ell(w)}(\Fl(d);\;\Zz)$ are a $\Zz$-basis for 
$R_n=H^*(\Fl(d);\:\Zz)$.

Alternatively, $R_n$ may be described in terms of Chern classes.  Let $U_i$ be the rank-$i$
vector bundle on $\Fl(n)$ whose fiber at a flag $E_{\bullet}$ is $E_i$, and let
$\xi_i=-c_1(U_{i}/U_{i-1})$.  Then $R_n$ is the quotient of $\Zz[\xi_1,\dots,\xi_n]$ by the
ideal generated by the elementary symmetric functions in the $\xi_i$.  The \defterm{Schubert
polynomials} $\Sch_w$ express the cohomology classes $[\Omega_w]$ of the Schubert varieties as
polynomials in the $\xi_i$.
Since $\Fl(n)$ is a smooth variety, its cohomology ring is the same as its Chow 
or intersection ring, allowing one to solve problems in enumerative geometry by means of
computations in $R_n$ involving Schubert polynomials.  This is the classical Schubert 
calculus.

The Schubert polynomials may be calculated using Demazure's divided difference operators;
see~\cite[pp.~170--173]{Ful97}.  We note here some special cases for later use, omitting the
calculations.  We write a permutation $w \in S_n$ in one-line notation
as a sequence $(w(1),w(2),\dots,w(n))$ (sometimes omitting the commas for brevity):
\begin{itemize}
\item If $w$ is the identity permutation, then $\Sch_w=1$ (the fundamental class, since
$\Omega_w$ is the whole flag variety).
\item If $w=(n,n-1,\dots,2,1)$ is the unique permutation of maximal length, then
$\Sch_w = \xi_1^{n-1}\xi_2^{n-2}\dots\xi_{n-1}$.  This is the cohomology class of a point
in $\Fl(n)$.
\item If $w$ is the transposition of $i$ with $i+1$, then $\Sch_{w} = 
\xi_1+\dots+\xi_i$.
\item If $w = (n,1,2,\dots,n-1)$, then $\Sch_{w} = \xi_1^d$.
\end{itemize}

The Schubert calculus may be extended to the partial flag manifold
    $$
    \Fl^{1,2}(d+1) \ = \ \left\{ F_1 \subset F_2 \subset \Cc^{d+1}
    \st \dim F_i = i \right\},
    $$
which is isomorphic to the picture space $\PSD(L_1)$.
The natural surjection $\Fl(d+1) \to
\Fl^{1,2}(d+1)$, forgetting the data $F_3,\dots,F_d$, induces a decomposition of 
$\Fl^{1,2}(d)$
into Schubert cells $\Omega_w^{\circ}$ of a form analogous to (\ref{schubert-decomp}), where 
$w \in S_{d+1}$ is a permutation such that
    \begin{equation} \label{relevant-permutations}
    w_3 > w_4 > \dots > w_{d+1}.
    \end{equation}
For such permutations, the Schubert polynomial $\Sch_{w} = [\Omega_w]$
involves only the variables $\xi_1$ and $\xi_2$.

More generally, if $G=(V,E)$ is an orchard, $v \in V$, and $e \in E(v)$,
then there is a fibration
    $$
    \pi = \pi_{v,e} : \ \PSD(G) \to \PSD(L_1) \isom \Fl^{1,2}(d+1)
    $$
sending a picture $\Pic$ to the partial flag $\Pic(v) \subset \Pic(e)$.  Thus
$\PSD(G)$ is a disjoint union of ``Schubert cells'' of the form
    $$
    \Omega_{\mathbf{w}}^{\circ} \ = \ \bigcap_{e \in E(v)} 
    \pi^{-1}(\Omega^{\circ}_{w_{v,e}})
    $$
indexed by tuples $\mathbf{w}$ of permutations $w_{v,e} \in S_{d+1}$. It is easy to
prove (by induction on $|E|$) that each $\Omega_{\mathbf{w}}^{\circ}$ is isomorphic to an
affine space, and to obtain conditions on the permutations $w_{v,e}$ for
$\Omega_{\mathbf{w}}^{\circ}$ to be nonempty.  We expect that in general
the ``orchard Schubert
variety'' $\Omega_{\mathbf{w}} = \overline{\Omega_{\mathbf{w}}^{\circ}}$ should be a union of
cells;  however, it is not clear how to describe the partial order analogous
to~(\ref{cell-closure}).

\begin{example} \label{quasi-bruhat-example}
Let $G=K_2$, with vertices $v_1,v_2$ and edge $e$.
For $i \in \{1,2\}$, there is a surjection $\pi_i: \PS^2(K_2) \to \Fl(3)$,
sending a picture $\Pic$ to the complete flag
$0 \subset \Pic(v_i) \subset \Pic(e) \subset \Cc^3$.
Thus $\PS^{2}(K_2)$ decomposes into Schubert cells
    $$
    \Omega_\sigma^{\circ} \ = \ \pi_{1}^{-1}(\Omega_{\sigma_1}) \;\cap\;
    \pi_{2}^{-1}(\Omega_{\sigma_2})
    $$
where $\sigma = (\sigma_1,\sigma_2) \in S_3 \x S_3$. In fact, $\Omega_\sigma^{\circ}$ is nonempty
if and only if $\sigma_1(3)=\sigma_2(3)$.  One may verify that the closure of a cell is indeed
a union of cells, and that the closure order analogous to~(\ref{cell-closure}) is
given by the following diagram:
    $$
    
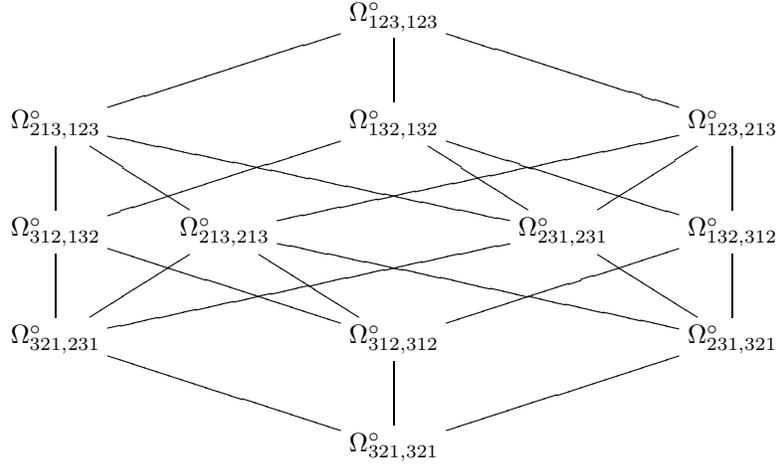
\xymatrix{
    && \Omega^\circ_{123,123} \ar@{-}[dll] \ar@{-}[d] \ar@{-}[drr] \\
    \Omega^\circ_{213,123} \ar@{-}[d] \ar@{-}[dr] \ar@{-}[drrr]
    && \Omega^\circ_{132,132} \ar@{-}[dll] \ar@{-}[dr] \ar@{-}[drr]
    && \Omega^\circ_{123,213} \ar@{-}[dlll] \ar@{-}[dl] \ar@{-}[d] \\
    \Omega^\circ_{312,132} \ar@{-}[d] \ar@{-}[drr]
    & \Omega^\circ_{213,213} \ar@{-}[dl] \ar@{-}[dr] \ar@{-}[drrr]
    && \Omega^\circ_{231,231} \ar@{-}[dlll] \ar@{-}[dr]
    & \Omega^\circ_{132,312} \ar@{-}[dll] \ar@{-}[d] \\
    \Omega^\circ_{321,231} \ar@{-}[drr]
    && \Omega^\circ_{312,312} \ar@{-}[d]
    && \Omega^\circ_{231,321} \ar@{-}[dll] \\
    && \Omega^\circ_{321,321} } 
    $$
Note that this is strictly weaker than the product of two copies of the strong Bruhat
order.  For instance, $231>213$ in Bruhat order, but the cells $\Omega_{231,231}^{\circ}$ 
and $\Omega_{213,213}^{\circ}$ both have complex dimension~$2$, hence are incomparable.
\end{example}

We can use the cell decomposition of $\PSD(G)$ to extend the methods of Schubert calculus to
orchards, using the presentation of the cohomology ring given in Theorem~\ref{orchard-thm}.
That is, we can count the number of pictures of an
orchard~$G=(V,E)$ meeting a given list of hyperplanes in $\Pp^d$ in certain ways.
For $v \in V$ and $e \in E(v)$, the map $\pi=\pi_{v,e}: \PSD(G) \to \Fl^{1,2}(d+1)$ induces a
graded ring homomorphism
    $$
    \pi^* : \ H^*(\Fl^{1,2}(d+1)) \to H^*(\PSD(G))
    $$
with the property that
    \begin{equation} \label{pullback-subvariety}
    \pi^*[Z] \ = \ [\pi^{-1}Z]
    \end{equation}
for all $Z \subset \Fl^{1,2}(d+1)$.  In particular, suppose that $w \in S_{d+1}$
is a permutation satisfying~(\ref{relevant-permutations}).
Then we can calculate the cohomology class
of $Y = \pi^{-1}(\Omega_{w})$ by evaluating the Schubert
polynomial $\Sch_w = \Sch_w(\xi_1,\xi_2)$ at $\xi_1=x_v$, \ $\xi_2=z_e-x_v$.  For example,
nine of the twelve varieties $\Omega_{\sigma_1,\sigma_2} \subset \PS^2(K_2)$ are pullbacks of
Schubert varieties in $\Fl(3)$ under $\pi_1$ or $\pi_2$ (the exceptions are 
$\Omega_{123,123}$,
$\Omega_{132,132}$, and $\Omega_{231,231}$).

Before giving an example of how this theory may be applied
to solve a problem in enumerative geometry,
we need one final ingredient---the cohomology class of a point.

\begin{proposition} \label{pointclass}
Let $G=(V,E)$ be an orchard.  With respect to the presentation~(\ref{concise-coho}), the 
cohomology class of a point in $\PSD(G,E)$ is
    $$
    \prod_{v \in V} \left( x_v^d \prod_{\text{loops}~e \in E(v)} h_{d-1}(x_v,y_e-x_v) \right).
    $$
\end{proposition}

\begin{proof}
Suppose first that $G$ is a forest (that is, it has no loops).
Let $v \in V$ and $p \in \Pp^d$.  Define
    $$
    Z_v \ = \ \{ \Pic \in \PSD(G) \st \Pic(v)=p\}.
    $$
Let $\pi : \PSD(G) \to \Pp^d$ be the map
sending $\Pic$ to $\Pic(v)$.  Then
    $$
    [Z_v] \ = \ [\pi^{-1}(p_0)] \ = \ \pi^*[p_0].
    $$
The cohomology ring of $\Pp^d$ is $\Zz[\xi]/\ideal{\xi^d}$,
where $\xi=c_1(\EE^*)$ and $\EE$ is 
the tautological line bundle on $\Pp^d$.  Then $\pi^{-1}\EE=\LL_{v}$, so
by~(\ref{chern-pullback}) and the previous equation we have $[Z_v]=x_v^d$.  The
discrete cellule is dense in $\PSD(G)$, because trees are $d$-parallel independent
for all $d$ by~\cite[Theorem 8.2.2]{Whiteley}.
This amounts to saying that the subvarieties $Z_v$ of
$\PSD(G)$ intersect transversely in a point.  Hence the cohomology class of a point is
$\prod_{v \in V} x_v^d$.

Now, suppose that $G$ contains one or more loops.  Let $q \in \Pp^d \sm \{p\}$.  For each 
loop $e \in E(v)$, define
    $$
    Y_e \ = \ \{\Pic \in \PSD(G) \st q \in \Pic(e)\}.
    $$
Then $Y_e = {\pi_{v,e}}^{-1}(\Omega_w)$, where $w = (d+1,1,d,d-1,\dots,2) \in
S_{d+1}$.  Using the special cases of Schubert polynomials mentioned previously and the
Demazure divided difference operators (we omit the details), one may show
that $[Y_e] = \Sch_{w} = h_{d-1}(x_v,y_e-x_v)$.  Then the collection of subvarieties
    $$
    \{Z_v \st v \in V \} \ \cup \ \{Y_e \st e \in E \text{ is a loop}\}
    $$
intersects transversely in a point, which implies the desired result.
\end{proof}

\begin{example} \label{coho-calc-example}
Let $G$ be the tree with vertices $V=\{1,2,3\}$ and edges $E=\{12,13\}$:
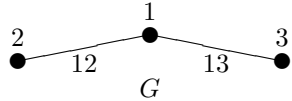
    \begin{center}
    \begin{picture}(100,50)
    \put(  0,20){\circle*{6}}   \put(  0,29){\makebox(0,0){$2$}}
    \put( 50,30){\circle*{6}}   \put( 50,39){\makebox(0,0){$1$}}
    \put(100,20){\circle*{6}}   \put(100,29){\makebox(0,0){$3$}}
    \put(  0,20){\line(5, 1){50}}   \put( 25,19){\makebox(0,0){$12$}}
    \put( 50,30){\line(5,-1){50}}   \put( 75,19){\makebox(0,0){$13$}}
    \put( 50,10){\makebox(0,0){$G$}}
    \end{picture}
    \end{center}
Let $A_1,A_2,A_3 \subset \Pp^3$ be planes, and let
$A_4, \dots, A_9 \subset \Pp^3$ be lines, with the
collection $\{A_i\}$ in general position.  We will calculate the number of 
pictures of $G$ in $\Pp^3$ satisfying the conditions
    \begin{equation} \label{schubert-conditions}
    \begin{aligned}
    \Pic(i) \in A_i \quad & \text{for } i=1,2,3, \\
    \Pic(12) \cap A_{i} \neq \emptyset \quad & \text{for } i=4,5,6, \\
    \Pic(13) \cap A_{i} \neq \emptyset \quad & \text{for } i=7,8,9.
    \end{aligned}
    \end{equation}

For $i=1,\dots,9$, let $Y_i$ be the subvariety of $\PS^3(G)$ consisting of pictures $\Pic$ for
which the condition involving $A_i$ is satisfied.  Then the problem is to determine the
cardinality of $Y = \bigcap_i Y_i$. Each $Y_i$ is the pullback of some $\Omega_w \subset
\Fl^{1,2}(\Cc^4)$, so its cohomology class is a Schubert polynomial in the variables
$x_1,x_2,x_3,z_{12},z_{13}$.  For instance,
    \begin{eqnarray*}
    \left[Y_1\right] &=& \left[{\pi_{1,12}}^{-1}(\Omega_{2134})\right] \ = \ 
\Sch_{2134}(x_1,z_{12}-x_1)
    \ = \ x_1 \qquad \text{and} \\
    \left[Y_4\right] &=& \left[{\pi_{1,12}}^{-1}(\Omega_{1324})\right] \ = \
    \Sch_{1324}(x_1,z_{12}-x_1) \;=\; z_{12}.
    \end{eqnarray*}
The other classes $[Y_i]$ may be calculated similarly.  In summary,
    \begin{align*}
    [Y_1] &= x_1, & [Y_4]=[Y_5]=[Y_6] &= z_{12}, \\
    [Y_2] &= x_2, & [Y_7]=[Y_8]=[Y_9] &= z_{13}, \\
    [Y_3] &= x_3.
    \end{align*}
Therefore $[Y] = x_1x_2x_3z_{12}^3z_{13}^3$.
By Proposition~\ref{pointclass}, the cohomology class of a point in 
$\PS^{3}(G)$ is $(x_1x_2x_3)^3$.  Since
$x_1x_2x_3z_{12}^3z_{13}^3 = 4(x_1x_2x_3)^3$
in $H^*(\PS^{3}(G))$, we conclude that $|Y|=4$.  That is,
there are four pictures of~$G$ satisfying 
the conditions~(\ref{schubert-conditions}).
(For this and many similar computations, we used the 
computer algebra system {\it Macaulay}~\cite{Macaulay}).
\end{example}

This calculation can be explained purely geometrically, in the spirit of the classical
Schubert calculus.  First, we specialize to the case that the lines $A_4$ and $A_5$ meet in a
point, as do $A_7$ and $A_8$. It is a fact that if there are four solutions to the
constraints~(\ref{schubert-conditions}) under this specialization, then there are four
solutions without it.  (In Schubert's terminology, this is the ``principle of conservation of
number''; see~\cite{KL72}.)

We will describe the set $J_{12}$ of locations for $\Pic(1)$ for which the
conditions on $\Pic(12)$ can be satisfied.  That is, $J_{12}$ consists of all points $p \in
A_1$ such that there exists a line $\ell$ through $p$ meeting each of $A_4$, $A_5$, $A_6$
nontrivially. There are two possibilities:

{\bf Case 1:} $p \in A_1 \cap Q$, where $Q$ is the plane containing $A_4$ and $A_5$.  Then
there is precisely one possibility for the line $\ell$: it must be the unique line determined
by $p$ and the point $A_6 \cap Q$.  Note that $\ell \subset Q$, so both $\ell \cap A_4$ and
$\ell \cap A_5$ are nonempty.

{\bf Case 2:} $p \in A_1 \cap R$, where $R$ is the plane determined by the point $a = A_4 \cap 
A_5$ and the line $A_6$.  Again, this choice of $p$ determines $\ell$ uniquely: it is the line 
determined by $p$ and $a$.  Note that $\ell \subset R$, so $\ell \cap A_6 \neq \emptyset$.

Thus $J_{12}$ is the union of two lines in $A_1$.  By an identical argument, the set $J_{13}$
of points for which the conditions on $\Pic(13)$ can be satisfied is also the union of two
lines in $A_1$.  Therefore $J = J_{12} \cap J_{13}$ consists of four points in $A_1$.  For
each point $p \in J$, there is a unique picture of~$G$ satisfying all the conditions
of~(\ref{schubert-conditions}), with $\Pic(1)=p$.  The case analysis above tells us how to
choose $\Pic(12)$ and $\Pic(13)$, and $\Pic(2)$ and $\Pic(3)$ must be respectively $\Pic(12)
\cap A_2$ and $\Pic(13) \cap A_3$.

This geometric result verifies \textit{ex post facto} that the subvarieties
$Y_i$ of Example~\ref{coho-calc-example} intersect transversely, so that
the above cohomological calculation is valid.  It would be interesting to discover
whether such transversality holds for all orchard Schubert varieties.


\section{Some Open Problems}

1. What aspects of the graph- or matroid-theoretic structure
of $G$ (other than its $d$-parallel behavior)
can be read off from the Poincar\'e series of $\PSD(G)$?
For instance, can the Tutte polynomial itself be recovered from
the Poincar\'e series?  This seems intuitively unlikely because
$\PSD(G)$ involves one fewer variable than $\Tutte_G(x,y)$.
However, the author's
experimentation has thus far produced no counterexample.

2. Theorem~\ref{d-parallel-thm}
may be read as saying that the $d$-parallel behavior of a graph is encapsulated
in the structure of the corresponding graphic matroid (which is less
information than the structure of the graph itself!)  Accordingly,
given an arbitrary matroid $M$, can one define the ``$d$-parallel matroid''
of $M$, with the same ground set, purely in terms of the Tutte polynomial?
(This idea was suggested to the author independently by M.~Haiman and V.~Reiner.) 
If so, what is the geometric meaning of such a combinatorial object?

3. Some of the material of Section~\ref{orchard-section} may warrant further
investigation.  Most glaring is the lack of a Schubert calculus for graphs
other than orchards, for which the picture space is singular.  In addition,
one might study the ``quasi-Bruhat'' poset associated with an orchard, as in
Example~\ref{quasi-bruhat-example} above.


\bibliographystyle{amsplain}
\bibliography{biblio}

\end{document}